\documentclass[english]{amsart}
\usepackage{babel}
\usepackage{caption}
\usepackage{caption}
\usepackage{subcaption}
\usepackage{tikz}

\usepackage{hyperref}

\newcommand{\kom}[1]{}
\renewcommand{\kom}[1]{{\bf [#1]}}

\def\vint_#1{\mathchoice%
          {\mathop{\kern 0.2em\vrule width 0.6em height 0.69678ex depth -0.58065ex
                  \kern -0.8em \intop}\nolimits_{\kern -0.4em#1}}%
          {\mathop{\kern 0.1em\vrule width 0.5em height 0.69678ex depth -0.60387ex
                  \kern -0.6em \intop}\nolimits_{#1}}%
          {\mathop{\kern 0.1em\vrule width 0.5em height 0.69678ex
              depth -0.60387ex
                  \kern -0.6em \intop}\nolimits_{#1}}%
          {\mathop{\kern 0.1em\vrule width 0.5em height 0.69678ex depth -0.60387ex
                  \kern -0.6em \intop}\nolimits_{#1}}}
                  
                  \newcommand{\aveint}[2]{\mathchoice%
          {\mathop{\kern 0.2em\vrule width 0.6em height 0.69678ex depth -0.58065ex
                  \kern -0.8em \intop}\nolimits_{\kern -0.45em#1}^{#2}}%
          {\mathop{\kern 0.1em\vrule width 0.5em height 0.69678ex depth -0.60387ex
                  \kern -0.6em \intop}\nolimits_{#1}^{#2}}%
          {\mathop{\kern 0.1em\vrule width 0.5em height 0.69678ex depth -0.60387ex
                  \kern -0.6em \intop}\nolimits_{#1}^{#2}}%
          {\mathop{\kern 0.1em\vrule width 0.5em height 0.69678ex depth -0.60387ex
                  \kern -0.6em \intop}\nolimits_{#1}^{#2}}}

\usepackage{color}
\usepackage{epsfig}

\newcommand{\R}{\mathbb{R}}

\def\1{\raisebox{2pt}{\rm{$\chi$}}}

\newcommand{\abs}[1]{\left| #1 \right|}

\newcommand{\ol}{\overline}
\newcommand{\Om}{\Omega}
\newcommand{\I}{\textrm{I}}
\newcommand{\II}{\textrm{II}}

\renewcommand{\P}{\mathbb{P}}

\newcommand{\vp}{\varphi}

\newcommand{\kint}{\vint}

\newcommand{\eps}{\varepsilon}

\theoremstyle{plain}
\newtheorem{definition}{Definition}[section]

\newtheorem{theorem}[definition]{Theorem}
\newtheorem*{theorem*}{Main Theorem}
\newtheorem{corollary}[definition]{Corollary}

\newtheorem{remark}[definition]{Remark}

\theoremstyle{definition}

\theoremstyle{remark}

\numberwithin{equation}{section}

\begin{document}

\title[]
{\bf Game-theoretic approach to H\"{o}lder regularity for PDEs involving eigenvalues of the Hessian}

\author[Blanc] {Pablo Blanc}
\address{Department of Mathematics and Statistics, University of Jyv\"{a}skyl\"{a}, P.O. Box 35, FI-40014 Jyv\"{a}skyl\"{a}, Finland}
\email{pblanc@dm.uba.ar}

\author[Han]{Jeongmin Han}
\email{jehan@jyu.fi}

\author[Parviainen]{Mikko Parviainen}
\email{mikko.j.parviainen@jyu.fi}

\author[Ruosteenoja]{Eero Ruosteenoja}
\email{eero.ruosteenoja@jyu.fi}

\date{\today}

\maketitle

\begin{abstract}
We prove a local H\"{o}lder estimate for any  exponent $0<\delta<\frac 12$ for solutions of the dynamic programming principle \begin{align*}u^\eps (x) =\sum_{j=1}^n \alpha_j\inf_{\dim(S)=j}\sup_{\substack{v\in S\\ |v|=1}}
\frac{ u^\varepsilon (x + \eps v) + u^\eps (x - \eps v)}{2}
\end{align*} 
with $\alpha_{1},\alpha_{n}>0$
and $\alpha_{2},\cdots,\alpha_{n-1}\ge0$.  
The proof is based on a new coupling idea from game theory. As an application, we get the same regularity estimate for viscosity solutions of the PDE $$\sum_{i=1}^n \alpha_i\lambda_i(D^2u)=0,$$ where $\lambda_1(D^2 u)\leq\cdots\leq \lambda_n(D^2 u)$ are the eigenvalues of the Hessian.

\end{abstract}

\section{Introduction}

\subsection{Main results} In this paper, we show local H\"{o}lder regularity for solutions of the following \textit{Dynamic Programming Principle} (DPP)
\begin{align} \label{dpp2}
u^\eps (x) =\sum_{j=1}^n \alpha_j\inf_{\dim(S)=j}\sup_{\substack{v\in S\\|v|=1}}
\frac{ u^\eps (x + \eps v) + u^\eps (x - \eps v)}{2}
\end{align}
in a bounded domain $\Omega \subset \mathbb{R}^{n}$,
where $\alpha_1,...,\alpha_n \geq 0$, $\min(\alpha_1,\alpha_n)>0$, and $\sum_{j=1}^n \alpha_j=1$. 
Our main result is the following, restated as Theorem \ref{thm:comb-holder}.

\begin{theorem*}
Let $u^{\eps}$ be a function satisfying the DPP \eqref{dpp2} in a bounded domain $\Omega$.
Then for any $0<\delta<\frac{1}{2}$ and $x,z\in B_{r}$ with $B_{2r} \subset\Omega$, there exists a constant $C=C(\delta, \alpha_{1}, \alpha_{n})>0$ such that
\begin{align}\label{mainest}|u^{\eps}(x)-u^{\eps}(z)| \le C||u^{\eps} ||_{L^{\infty}(B_{2r})}\bigg(\frac{|x-z|^{\delta}}{r^{\delta}}+\frac{\eps^{\delta}}{r^{\delta}}\bigg). 
\end{align} 
\end{theorem*}

That the above theorem holds for any $0<\delta<\frac{1}{2}$ is explicitly obtained in the proof in (\ref{eq:choice-of-delta}).

The DPP \eqref{dpp2} has a connection to a certain  PDE involving eigenvalues of the Hessian. Indeed, under certain regularity assumptions for the boundary of the domain, when $\eps \rightarrow 0$, solutions of \eqref{dpp2} converge uniformly to the unique viscosity solution of the following PDE,
\begin{equation}
\label{eq:mainN}
\sum_{i=1}^n \alpha_i\lambda_i(D^2u)=0,
\end{equation}
where  $\lambda_1(X)\leq\cdots\leq \lambda_n(X)$ are the ordered eigenvalues of $X\in S(n)$, the set of $n\times n$ real symmetric matrices.

As a consequence of our main result, we obtain the same H\"{o}lder estimate for any $0<\delta<\frac{1}{2}$ for viscosity solutions of this PDE.
For the proof of the following corollary, 
see Section 2.3.

\begin{corollary}
\label{coro}
Let $u$ be the viscosity solution of \eqref{eq:mainN} in a bounded domain $\Omega$. Then for any $0<\delta<\frac{1}{2}$ and $x,z\in B_{r}$ with $B_{2r} \subset\Omega$, there exists a constant $C=C(\delta,\alpha_{1}, \alpha_{n})>0$ such that
$$|u (x)-u(z)| \le C||u||_{L^{\infty}(B_{2r})} \frac{|x-z|^{\delta}}{r^{\delta}}. $$ 
\end{corollary}

We remark that $\sum_{i=1}^n \alpha_i\lambda_i =0 $ satisfies Pucci type inequalities, and thus we can get H\"{o}lder regularity from the general theory directly, even though the exponent is not explicitly given (see the beginning of Section 4).
 
For this equation with lower order terms, Ferrari and Vitolo \cite{ferrari2020regularity} used methods from the viscosity theory to study ABP, Harnack and H\"{o}lder estimates, and later Ferrari and Galise \cite{ferrari2021regularity} showed $C^{0,\delta}$-regularity for $ \delta=1-\frac{\alpha_{1}+\alpha_{n}}{(\sqrt{\alpha_{1}}+\sqrt{\alpha_{n}})^{2}}\in (0, \frac 12 ]$. We note that $\delta= \frac{1}{2}$ holds if and only if $\alpha_{1}=\alpha_{n}$. 
Also observe that $\min(\alpha_1,\alpha_n)>0$ is necessary in \cite{ferrari2020regularity} as well as in our main result.

 The DPPs of type (\ref{dpp2}) model competition of two players, and one can relate these games to different applications. For example, in the context of related tug-of-war games, they have been suggested in connection to the option pricing problem with market manipulation \cite{nystromp17}. 
 To be more precise,  for a given boundary data, the solution of the DPP \eqref{dpp2} is also the value function of a two-player zero-sum stochastic game, the rules of which can be read from the DPP. We will describe the game in more detail in the next section, but informally, a token is placed at $x\in \Omega$, $\alpha_j$ is the probability that the number $j$ is chosen.
Then the player aiming to minimize the value chooses a $j$-dimensional subspace of $\R^n$, and finally, the player aiming to maximize the value chooses a unit  vector from that subspace. Then the token is moved an $\eps$-step either to the direction or the opposite direction of the vector, with equal probabilities. The game continues until the token is moved outside of $\Omega$, and the player choosing subspaces pays the amount given by the boundary payoff function to the other player.

To the best of our knowledge, there are no prior works studying local regularity of DPPs or games related to fully nonlinear PDEs involving eigenvalues of the Hessian.

The game that we just described is connected to the PDE \eqref{eq:mainN}.
This connection has been studied in detail by Blanc and Rossi \cite{blancr19} for the PDE $\lambda_j(D^2 u)=0$, where $j\in \{1,...,n\}$.
See also \cite{brustadlm20} for a game associated to the Dominative $ p $-Laplacian and \cite{hoeg2020control,blanc2020evolution} for games associated to parabolic versions of these equations.

The rest of the paper is organized as follows. In the following subsection, we give a more detailed idea of our proof method. In Section 2 we give some preliminary definitions and results for viscosity solutions. In Section 3 we prove the main result for the special case $\alpha_1= \alpha_n=\frac 12 $. In Section 4 we prove the main theorem.

\subsection{Method of the proof}
Although the ideas behind the proof of our main theorem stem from games, we do not use methods from stochastic game theory. Instead, our starting point is the coupling method introduced by Luiro and Parviainen \cite{luirop18} in the context of tug-of-war games. However, a direct application of their method does not seem to work in our case, so we need a new type of coupling, which is the main novelty of this work.

To give an idea of the proof, for simplicity we will discuss a special case $\alpha_1= \alpha_n = \frac 12 $, in which case the DPP \eqref{dpp2} can be written as 
\begin{align*}
u^\eps (x) =
\frac{1}{2}\sup_{|v|=1} \left\{ \frac{ u^\eps (x + \eps v) + u^\eps (x - \eps v)}{2}\right\}+
\frac{1}{2}\inf_{|w|=1}\left\{ \frac{ u^\eps (x + \eps w) + u^\eps (x - \eps w)}{2}\right\}.
\end{align*}
The starting point of the coupling is to define a $2n$-dimensional game related to the DPP. Notice that the function 
$g:\, \Omega\times \Omega\rightarrow \R$ given by \[ g(x,z)=u^\eps(x)-u^\eps(z)\] can be written as a solution of a suitable DPP in $\R^{2n}$ as follows,
\[
\begin{split}
& g(x,z)=u^\eps(x)-u^\eps(z)\\
&=
\frac{1}{2}\sup_{|v|=1}\left\{ \frac{ u^\eps (x + \eps v) + u^\eps (x - \eps v)}{2}\right\}
+
\frac{1}{2}\inf_{|\tilde v|=1}\left\{ \frac{ u^\eps (x + \eps\tilde v) + u^\eps (x - \eps\tilde v)}{2}\right\}\\
&\quad -
\frac{1}{2}\sup_{|\tilde w|=1}\left\{ \frac{ u^\eps (z + \eps \tilde w) + u^\eps (z - \eps \tilde w)}{2}\right\}
-
\frac{1}{2}\inf_{|w|=1}\left\{ \frac{ u^\eps (z + \eps w) + u^\eps (z - \eps w)}{2}\right\}\\
&=
\frac{1}{2}\sup_{\substack{|v|=1\\|w|=1}}\left\{ \frac{ u^\eps (x + \eps v) + u^\eps (x - \eps v)-u^\eps (z + \eps w) - u^\eps (z - \eps w)}{2}\right\}\\
&\quad +
\frac{1}{2}\inf_{\substack{|\tilde v|=1\\|\tilde w|=1}}\left\{ \frac{ u^\eps (x + \eps\tilde v) + u^\eps (x - \eps\tilde v)-u^\eps (z + \eps \tilde w) - u^\eps (z - \eps \tilde w)}{2}\right\}.\\
\end{split}
\]
The potential of this DPP is to transform the question of regularity of $u^{\eps}$ to the question of the absolute size of $g$. The heuristic idea is to introduce a suitable stochastic game in $\Omega\times \Omega$, where we aim to move the two tokens to the diagonal set $$T:=\{(x,z)\in \Omega\times \Omega\, :\, x=z\}, $$ where $g=0$, before our opponent can move the tokens outside of the set $\Omega\times \Omega$.
Observe that 
\begin{align*}
& \frac{ u^\eps (x + \eps v) + u^\eps (x - \eps v)-u^\eps (z + \eps w) - u^\eps (z - \eps w)}{2}\\
& \quad =\frac{ g((x,z) + \eps (v,w)) +g((x,z) - \eps (v,w)) }{2}.
\end{align*}
Following the idea of Luiro and Parviainen, we could consider a $2n$-dimensional game where each player (both with probability $\frac12$) gets to choose $v$ and $w$ and then the tokens move to $(x,z) + \eps (v,w)$ or to $(x,z) - \eps (v,w)$, each possibility with probability one half. If we set the boundary values of our game to be 0 on $T$ and $2\sup u^\eps$ in $\R^{2n}\setminus (\Omega\times \Omega)$ and could prove that $|g(x,z)|\leq C|x-z|^\delta$, we would get a desired H\"older estimate for the function $u^\eps$.

\begin{figure}
    \centering
    \begin{subfigure}[b]{0.5\textwidth}
    	\centering
    	\renewcommand{\thesubfigure}{A}
    	\begin{tikzpicture}[scale=0.5]
    	\node at (0,0,0)[circle,fill,inner sep=1pt,label = left:$x$]{};
    	\draw [-] (0,0,0) -- (4,0,0);
    	\node at (4,0,0)[circle,fill,inner sep=1pt,label = right:$y$]{};
    	\draw [->] (0,0,0) -- (0,2,0) node [midway, left] {$v$};
    	\draw [->] (4,0,0) -- (4,-2,0) node [midway, right] {$w$};
   		\end{tikzpicture}
        \caption{A case where the rules work out.}
        \label{fig:a}
    \end{subfigure}
    ~ 
    \begin{subfigure}[b]{0.4\textwidth}
    	\centering
    	\renewcommand{\thesubfigure}{B}
		\begin{tikzpicture}[scale=0.5]
    	\node at (0,0,0)[circle,fill,inner sep=1pt,label = left:$x$]{};
    	\draw [-] (0,0,0) -- (4,0,0);
    	\node at (4,0,0)[circle,fill,inner sep=1pt,label = right:$y$]{};
    	\draw [->] (0,0,0) -- (0,2,0) node [midway, left] {$v$};
    	\draw [->] (4,0,0) -- (4,0,2) node [right, below] {$w$};
   		\end{tikzpicture}
        \caption{A bad case.}
        \label{fig:b}
    \end{subfigure}
    \caption{Two choices of vectors when the tokens are at $(x,z)$.}
    \label{fig}
\end{figure}
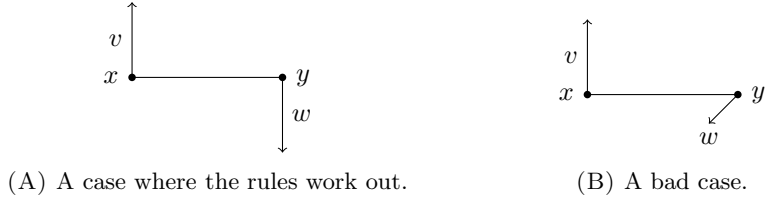

Unfortunately, these rules for the $2n$-dimensional game do not seem to be suitable to obtain regularity estimates in our case. The problem is that our opponent can force the tokens away by choosing $w=-v$ normal to $x-z$.
Observe that if the new position of the tokens is given by $(\tilde x,\tilde z)=(x,y)+\eps(v,w)$, we get $|\tilde x-\tilde z|^2= |x-z|^2+4\eps^2$.
The same holds for $(\tilde x,\tilde z)=(x,y)+\eps(-v,-w)$.

Observe that it also holds 
\begin{align*}
& \frac{ u^\eps (x + \eps v) + u^\eps (x - \eps v)-u^\eps (z + \eps w) - u^\eps (z - \eps w)}{2}\\
& \quad =\frac{ g((x,z) + \eps (v,-w)) +g((x,z) - \eps (v,-w)) }{2}.
\end{align*}
Again the rules that follow from this formula allow our opponent to force the tokens away from each other, in this case by choosing $w=v$ normal to $x-z$.

In conclusion, with the rules that we have described, we do not get a suitable coupling. Our new idea is to couple the moves in one way or the other depending on the choice of the vectors $v$ and $w$. When the tokens are placed at $x$ and $z$, given $v$ and $w$ we define two rules moving the token:
\begin{enumerate}
\item
the token moves to $(x,z) + \eps (v,w)$ or to $(x,z) - \eps (v,w)$, each possibility with probability one half.
\item
the token moves to $(x,z) + \eps (v,-w)$ or to $(x,z) - \eps (v,-w)$, each possibility with probability one half.
\end{enumerate}

Let us define the $2n$-dimensional game. We toss a coin and the winner of the toss chooses two unitary vectors $v$ and $w$. Set $y=x-z$. We also write $v_{y^{\perp}}=v-\frac{\langle v, y \rangle}{\langle y,y \rangle}y$ and 
$w_{y^{\perp}}=w-\frac{\langle w, y \rangle}{\langle y,y \rangle}y$.
If $|v_{y^{\perp}}|^2+|w_{y^{\perp}}|^2>1$ and $\langle v_{y^{\perp}},w_{y^{\perp}}\rangle<0$, then the tokens move according to rule (ii). In any other case the tokens move according rule (i).

Define 
\begin{align} \label{defcapf}\begin{split}&F(x,z,v,w,g) \\ &=
\left\{ \begin{array}{ll}
\frac{ g((x,z) + \eps (v,-w)) +g((x,z) - \eps (v,-w)) }{2} & \textrm{if $|v_{y^{\perp}}|^2+|w_{y^{\perp}}|^2>1$ and $\langle v_{y^{\perp}},w_{y^{\perp}}\rangle<0$,}\\
\frac{ g((x,z) + \eps (v,w)) +g((x,z) - \eps (v,w)) }{2}& \textrm{otherwise.}\\
\end{array} \right.
\end{split}
\end{align}
Then we obtain the DPP for our $2n$-dimensional game
\begin{align} \label{add3}
g(x,z)=
\frac{1}{2}\sup_{\substack{|v|=1\\|w|=1}}
F(x,z,v,w,g)
+
\frac{1}{2}\inf_{\substack{ | v|=1\\ | w|=1}}
F(x,z, v, w,g).
\end{align}

Observe that in the case of Figure~\ref{fig:a}, that is when a player selects $w=-v$ normal to $x-z$, we have to move the tokens accordingly to rule (ii).
We have  $x + \eps v-(z-\eps w))=x - \eps v-(z+\eps w))=x-z$, and therefore the distance between the tokens is preserved.

Now, we consider the case in Figure~\ref{fig:b}, that is, $v$ and $w$ are normal to $x-z$ and also to each other.
If the new position of the tokens is given by $(\tilde x,\tilde z)=(x,z)+\eps(v,w)$, we get $|\tilde x-\tilde z|^2= |x-z|^2+2\eps^2$.
We get the same if $(\tilde x,\tilde z)=(x,z)+\eps(v,-w)$.
Then, the tokens are forced away from each other independently of what rule we select.
But still, observe that this growth is smaller than the one we were getting in the case of Figure~\ref{fig:a} when applying rule (i), since
$|\tilde x-\tilde z|^2= |x-z|^2+4\eps^2 $ in that case.
We claim that our choice of when to apply (i) or (ii) reduces the ability of the players to push the tokens away from each other, and this is the key of the matter.

\section{Preliminaries}

In this section, we include some preliminary results concerning solutions to the equation and the game.

\subsection{Viscosity solutions}

Let us begin by recalling the definition of viscosity solutions for \eqref{eq:mainN}.
We denote by $USC(\overline{\Omega})$ (resp. $LSC(\overline{\Omega})$)
the set of all upper (resp. lower) semicontinuous functions defined on $\overline{\Omega}$.
\begin{definition}Let $u:\bar\Omega \to \mathbb{R}$ be a function.
\begin{itemize}
\item[(a)] We say that $u\in USC(\ol \Om)$ is a viscosity subsolution to \eqref{eq:mainN} if for all $x \in \Omega$ and $\phi \in C^{2}(\Omega)$ such that $u(x)=\phi(x)$ and
$u(y)<\phi(y)$ for $y \neq x$ we have 
\begin{align*} 
\sum_{i=1}^n \alpha_i\lambda_i (D^{2}\phi)\ge 0 .
\end{align*}
\item[(b)] We say that $u\in LSC(\ol \Om)$ is a viscosity supersolution to \eqref{eq:mainN} if for all $x \in \Omega$ and $\phi \in C^{2}(\Omega)$ such that $u(x)=\phi(x)$ and
$u(y)>\phi(y)$ for $y \neq x$ we have 
\begin{align*}
\sum_{i=1}^n \alpha_i\lambda_i (D^{2}\phi) \le 0 .
\end{align*}
\end{itemize}
If $u$ is continuous and satisfies both of \emph{(a)} and \emph{(b)}, we say that $u$ is a viscosity solution to \eqref{eq:mainN}.
\end{definition}

Next, we prove comparison and thus uniqueness to our operator. It would follow from \cite{birindelli2020existence}, but here we give a simple alternative proof for this particular operator.

\begin{remark}
If $M$ is a Hermitian matrix, it is diagonalizable with real eigenvalues and by the Min-max Theorem those eigenvalues verify
\[
\lambda_j(M)=\inf_{\dim(S)=j}\sup_{|v|=1} \langle M v, v \rangle
\]
for $j=1,\dots,N$, and we can use this identity for the Hessian matrix.
\end{remark}

\begin{theorem}
Let $u_{1}\in USC(\ol \Om)$ be a viscosity subsolution (\ref{eq:mainN}), and $u_{2}\in LSC(\ol \Om)$ a viscosity supersolution to (\ref{eq:mainN}). 
\[
\begin{split}
\text{If }u_{1}\le u_{2}\text{ on } \partial \Om \text{, then } u_{1}\le u_{2} \text{ on }\ol \Om.
\end{split}
\]
\end{theorem}
\begin{proof}
First, we make a counter proposition
\[
\begin{split}
\sup_{\Om}(u_{1}-u_{2})=:\theta>0.
\end{split}
\]
\sloppy
Then observe that the equation can be written as  $$\sum_{i=1}^n \alpha_i\lambda_i(D^2u_{1})=\sum_{i=1}^{n} \alpha_i\inf_{\dim(S)=j}\sup_{|v|=1}  \langle D^2  u_{1} (x) v, v \rangle .$$
Moreover,  $\tilde u_{1}(x)=u_{1}(x)+\delta\abs{x}^2$ is a subsolution to 
\begin{align*}
\label{eq:strict}
\sum_{i=1}^{n} \alpha_i\inf_{\dim(S)=j}\sup_{|v|=1}  \langle D^2 \tilde  u_{1} (x) v, v \rangle &=  \sum_{i=1}^{n} \alpha_i\inf_{\dim(S)=j}\sup_{|v|=1}  \langle D^2(u_{1} (x) +\delta\abs{x}^2)v, v \rangle \\ &=2C\delta.
\end{align*}
Let 
\[
\begin{split}
\Phi(x,y)&=\tilde u_{1}(x)-u_{2}(y)-\vp(x,y)\\
&:=\tilde u_{1}(x)-u_{2}(y)-\frac{1}{2\eps}\abs{x-y}^2,
\end{split}
\]
 and let  $(x_\eps,y_\eps)$ be the maximum point on $\ol \Om\times \ol \Om$. The points  $x_\eps, y_\eps$ are not at the boundary in a bounded domain $\Om$ for small enough $\delta$ by the standard theory.
Then by the theorem of sums \cite{crandallil92} we have
\[
\begin{split}
(D_x\vp(x_\eps,y_\eps),X)\in \ol J^{2,+}\tilde u_{1}(x_\eps),\quad (-D_y\vp(x_\eps,y_\eps),Y)\in \ol J^{2,-}u_{2}(y_\eps),
\end{split}
\]
and $X\le Y$. Furthermore, let $\eta>0$, $S_j$ be a $j$-dimensional subspace of $\mathbb{R}^{n}$ and $v_{j} \in S_{j}$ with $\abs{v_j}=1$ be such that $$\sup_{v\in S_j,|v|=1}  \langle Y v, v \rangle \le \inf_{\dim(S)=j}\sup_{|v|=1}  \langle Y v, v \rangle+\eta$$ and $$\eta+\langle  X v_j, v_j \rangle\ge \sup_{v\in S_j,|v|=1} \langle  X v, v \rangle .$$
Then we have 
\[
\begin{split}
2C\delta &\le \sum_{i=1}^{n}\alpha_{i}\Big(\inf_{\dim(S)=j}\sup_{|v|=1} \langle  X v, v \rangle-  \inf_{\dim(S)=j}\sup_{|v|=1}  \langle Y v, v \rangle\Big)\\
&\le  \sum_{i=1}^{n}\alpha_{i}\Big(\inf_{\dim(S)=j}\sup_{|v|=1} \langle  X v, v \rangle- \sup_{v\in S_j,|v|=1}  \langle Y v, v \rangle\Big)+\eta\\
&\le \sum_{i=1}^{n}\alpha_{i}\Big(\sup_{v\in S_j,|v|=1} \langle  X v, v \rangle- \sup_{v\in S_j,|v|=1}  \langle Y v, v \rangle\Big)+\eta\\
&\le \sum_{i=1}^{n}\alpha_{i}\Big(\langle  X v_j, v_j \rangle- \sup_{v\in S_j,|v|=1}  \langle Y v, v \rangle\Big)+2\eta\\
&\le  \sum_{i=1}^{n}\alpha_{i}(\langle  X v_j, v_j \rangle-  \langle Y v_j, v_j \rangle)+2\eta\\
&\le  \sum_{i=1}^{n}\alpha_{i}\langle (X-Y) v_j, v_j \rangle+2\eta\le 2\eta,
\end{split}
\]
which is a contradiction for small enough $\eta>0$. 
\end{proof}
It might also be instructive to think some special cases. For example, for the first eigenvalue equation  $\lambda_1 (D^2 u) (x) := \inf_{|v|=1} \langle D^2 u (x) v, v \rangle = 0$,
choosing $\langle  Y v_0, v_0 \rangle =\inf_{|v|=1} \langle  Y v, v \rangle$, the key computation above reads as
\[
\begin{split}
2\delta &\le \inf_{|v|=1} \langle  X v, v \rangle-  \inf_{|v|=1}  \langle Y v, v \rangle\\
&\le \inf_{|v|=1}\langle  X v, v \rangle-  \langle Y v_0, v_0 \rangle+\eta\\
&\le  \langle  X v_0, v_0 \rangle- \langle Y v_0, v_0 \rangle+\eta\\
&\le \langle (X-Y) v_0, v_0 \rangle+\eta\le \eta,
\end{split}
\]
a contradiction.

Observe that uniqueness immediately follows from the comparison principle for the viscosity solutions of the boundary value problem 
\begin{align*}
\begin{cases}
\sum_{i=1}^n \alpha_i\lambda_i(D^2u)=0 & \text{ in }\Om,\\
u=g & \text{ on } \partial \Om, 
\end{cases}
\end{align*}
with given continuous boundary values $g:\partial \Om\to \R$. 
Also, observe that if the domain is strictly convex, we have that a plane can act as a barrier.
Then, the solution obtained by Perron's method turns out to be continuous up to the boundary.
A weaker condition for the existence of continuous solutions for smooth domains can be found in \cite{harveyl09}.

\subsection{Games}

A game associated with the equation $\lambda_j(D^2 u)=0$ was introduced in \cite{blancr19}. 
Here we modify the game so that it is associated with equation \eqref{eq:mainN}.
Next, we give the precise formulation of the game.

It is a two-player zero-sum game. Fix a domain $\Omega \subset \mathbb{R}^N$, $\eps>0$ and a final payoff function 
$G :\mathbb{R}^N \setminus \Omega \mapsto \mathbb{R}$. 
The rules of the game are the following:
the game starts with a token at an initial position $x_0 \in \Omega$ and develops in several rounds.
At the beginning of each round $j\in\{1,\dots,n\}$ is chosen at random such that $\P(j=i)=\alpha_i$ for each $i=1,\dots,n$. 
With this given value, Player~I chooses a subspace $S$ of dimension $j$ and subsequently Player~II a unitary vector $v\in S$.
Then the position of the token is moved to $x\pm \eps v$ with equal probabilities. 
The game continues until the token leaves the domain.
At this time that we call $\tau$, Player~I pays $G(x_\tau)$ to Player~II.
When the two players fix their strategies $S_I$ and $S_{II}$, we can compute the expected outcome as
$$
\mathbb{E}_{S_\I, S_{\II}}^{x_0} [G (x_\tau)].
$$
Then the value of the game for any $x_0 \in \Omega$ is defined as
\[
u^\eps(x_0)=\sup_{S_\I}\inf_{S_{\II}}\,
\mathbb{E}_{S_{\I},S_\II}^{x_0}\left[G(x_\tau)\right]=\inf_{S_{\II}}\sup_{S_\I}\,
\mathbb{E}_{S_{\I},S_\II}^{x_0}\left[G(x_\tau)\right],
\]
and verifies the DPP \eqref{dpp2}, that is
\begin{align} \label{pde1} 
u^\eps (x) =\sum_{j=1}^n \alpha_j\inf_{dim(S)=j}\sup_{\substack{v\in S\\ |v|=1}}
\frac{ u^\eps (x + \eps v) + u^\eps (x - \eps v)}{2}
\end{align}
for $x \in \Omega$ and $u^\eps(x) = G(x)$ for $x \not\in \Omega$,
see \cite{blancr19book}. 
Intuitively, the rules of the game can be seen from the DPP.
When Player I, who aims to minimize the value, chooses a subspace, she knows that Player II aims to choose from that subspace a unitary vector maximizing the average `$\eps$-step value'.

\subsection{Application to the PDE \eqref{eq:mainN}}

We give a brief explanation of the relation between \eqref{dpp2} and \eqref{eq:mainN}, and how to prove Corollary \ref{coro}
using the result of our main theorem.

If we assume that the domain is strictly convex, we obtain that
\begin{equation}
\label{eq:lim}
u^\eps\to u
\end{equation}
uniformly as $\eps\to 0$ where $u$ is the unique solution to \eqref{eq:mainN}.
Observe that in \cite{blancr19} a condition over the boundary is given for each $j$.
This condition is used to prove that the game value is asymptotic continuous near the boundary. 
It is in this step that we use that the domain to be strictly convex.
Then the convergence is obtained following the usual path, see \cite{MR2451291,manfredipr12,blancr19book}.
We use the asymptotic version of Arzela-Ascoli to pass to the limit, and using the DPP combined with the definition of viscosity solutions allows us to deduce that the limit is the unique viscosity solution.  
Observe that a weaker condition on the domain may be enough depending on $\alpha_i$ but we are not going to address this question here.

The connection between the DPP and the PDE can be intuitively seen by recalling that 
\[
\lambda_j(M)=\inf_{\dim(S)=j}\sup_{|v|=1} \langle M v, v \rangle
\]
and observing that
\[
\frac{ u^\eps (x + \eps v) + u^\eps (x - \eps v)}{2}
\approx
\eps^2\langle D^2u(x) v, v \rangle.
\]

Assume that the estimate \eqref{mainest} is provided for any $\eps>0$.
We observe that since a ball is strictly convex we obtain the convergence \eqref{eq:lim} in there.
By passing to the limit as $\eps\to 0$ in the main theorem, we obtain a H\"{o}lder estimate for the limit function $u$.
That is what we stated as Corollary~\ref{coro}.

\section{Regularity for a DPP related to $\frac12 \lambda_1(D^2 u)+\frac12 \lambda_n(D^2 u)=0$}

We first focus on the case given by $\alpha_1=\alpha_n=\frac{1}{2}$. 
Note that in this case $\alpha_2=\cdots=\alpha_{n-1}=0$, since we assumed $\sum_{j=1}^n \alpha_j=1$.
Then the DPP \eqref{dpp2} is simplified to
\begin{align} \label{dpp1}
u^\eps (x) =
\frac{1}{2}\sup_{|v|=1}\left\{ \frac{ u^\eps (x + \eps v) + u^\eps (x - \eps v)}{2}\right\}
+
\frac{1}{2}\inf_{|w|=1}\left\{ \frac{ u^\eps (x + \eps w) + u^\eps (x - \eps w)}{2}\right\}.
\end{align}

The game starts with a token at an initial position $x_0 \in \Omega$. At every round, a fair coin is tossed and the winner of the toss chooses a vector $v\in \mathbb{R}^n$ with $|v|=1$. Then the position of the token is moved to either $x_0 +\eps v$ or $x_0 -\eps v$, with equal probabilities. 
The game ends when the token leaves the domain and we define the game value as before. 
Our game value satisfies the DPP \eqref{dpp1}
for $x \in \Omega$, and $u^\eps(x) = G(x)$ for $x \not\in \Omega$.

In this section, we will obtain the regularity result for solutions to the DPP \eqref{dpp1}. 
As we have mentioned, we employ the method introduced in \cite{luirop18}.
For the readers' convenience, we will provide a full proof.

\begin{theorem}
\label{thm:main}
Let $u^{\eps}$ be a function satisfying the DPP \eqref{dpp1} in a bounded domain $\Omega$.
Then for any $0<\delta<\frac{1}{2}$ and $x,z\in B_{r}$ with $B_{2r} \subset\Omega$, there exists a constant $C=C(\delta)>0$ such that
$$|u^{\eps}(x)-u^{\eps}(z)| \le C||u^{\eps} ||_{L^{\infty}(B_{2r})}\bigg(\frac{|x-z|^{\delta}}{r^{\delta}}+\frac{\eps^{\delta}}{r^{\delta}}\bigg). $$ 
\end{theorem}

\begin{proof}
By considering $\tilde u(x)=u(rx)$ we can assume that $r=1$.
Also, without loss of generality, we assume that 
\begin{align} \label{as1} \sup_{(x,z) \in B_{2}\times B_{2}} (u^{\eps}(x)-u^{\eps}(z)) = \sup_{(x,z) \in B_{2}\times B_{2}} g(x,z)\le 1
\end{align}
by a suitable renormalization.

To obtain the desired estimate for the function, we will use the comparison function $f=f_1-f_2$ with
\[
\begin{split}
f_1(x,z)=C|x-z|^{\delta}+|x+z|^2\,
\end{split}
\]
and
\begin{equation*} f_{2}(x, z) = \left\{ \begin{array}{ll}
C^{2(N-i)} \eps^ {\delta}  & \textrm{if $(x, z) \in A_{i}$,}\\
0 & \textrm{if $|x-z|>N \eps / 10 $,}
\end{array} \right.
\end{equation*}
where $C>1$,
$\delta\in (0,1)$, and
$$A_{i} = \{ (x , z) \in \mathbb R ^ {2n} : (i-1) \eps / 10 < |x-z| \le i \eps / 10 \} $$ for $ i = 0, 1, ... , N$, where $N$ is a sufficiently large number to be determined later.

The first term in $f_1$ will give us the desired regularity estimate, and the second term ensures that the estimate holds at  $(B_2\times B_2)\setminus (B_1\times B_1)$.
It is typical for the solutions of `$\eps$-DPPs' that they are discontinuous at the $\eps$-scale. That is why we need a correction function $f_2$, designed to handle the case where the distance $|x-z|\approx \eps$.

We first observe that 
\begin{equation} \label{as2}
\begin{split}
u^{\eps}(x)-u^{\eps}(z)-f(x,z)
&= u^{\eps}(x)-u^{\eps}(z)-f_1(x,z)+f_2(x,z)\\ 
&\leq \max f_2
\le C^{2N}\eps^{\delta}.
\end{split}
\end{equation}
for $(x,z) \in (B_{2}\times B_{2}) \backslash ((B_{1}\times B_{1}) \backslash T)$, where we are using that $f_1\geq 1$ in $(B_{2}\times B_{2}) \backslash (B_{1}\times B_{1})$.
If we prove that 
\begin{align} \label{int2}
\sup_{(x,z) \in  B_{2}\times B_{2}  } (u^{\eps}(x)-u^{\eps}(z)-f(x,z))\le C^{2N}\eps^{\delta},
\end{align}
the result follows as we can assume without loss of generality with a suitable translation that $x=-z$ 
since in this case we can obtain 
$$u^{\eps}(x)-u^{\eps}(-x) \le  2C^{2N}(|x|^{\delta}+\eps^{\delta}) $$
from
$f(x,-x)\le C(2|x|)^{\delta}.$

We assume, for the sake of contradiction,  that 
\begin{align} \label{as3}
M := \sup_{(x,z) \in ( B_{1}\times B_{1}) \backslash T} (u^{\eps}(x)-u^{\eps}(z)-f(x,z)) > C^{2N}\eps^{\delta}.
\end{align}
In this case, we have
\begin{align}  \label{int1}
M= \sup_{(x,z) \in  B_{2}\times B_{2}  } (u^{\eps}(x)-u^{\eps}(z)-f(x,z)) .
\end{align}
Consider an arbitrary small number $\eta >0$. Then we can choose $(x_{1},z_{1}) \in  (B_{1}\times B_{1}) \backslash T  $ such that
\begin{align*}
M \le u^{\eps}(x_{1})-u^{\eps}(z_{1})-f(x_{1},z_{1})+\eta.
\end{align*}

Recall \eqref{defcapf}. From \eqref{int1}, we observe that $$ u^{\eps}(x)-u^{\eps}(z) =g(x,z) \le M +f(x,z)$$  
for any $(x,z) \in B_{2}\times B_{2} $.
If two unit vectors $v$ and $w$ satisfy $|v_{y^{\perp}}|^2+|w_{y^{\perp}}|^2>1$ and $\langle v_{y^{\perp}},w_{y^{\perp}}\rangle<0$,
we have
\begin{align*}
 F(x_{1},z_{1},v,w,g) &
= \frac{ g((x_{1},z_{1}) + \eps (v,-w)) +g((x_{1},z_{1}) - \eps (v,-w)) }{2}
\\ & \le \frac{ \big( M+ f((x_{1},z_{1}) + \eps (v,-w))\big) + \big( M + f((x_{1},z_{1}) - \eps (v,-w)) \big) }{2}
\\ & = M + \frac{ f((x_{1},z_{1}) + \eps (v,-w)) +f((x_{1},z_{1}) - \eps (v,-w)) }{2}
\\& \le M + F(x_{1},z_{1},v,w,f)
\end{align*}
since $x_{1}\pm \eps v, z_{1}\pm \eps w $ are still in $ B_{2}$.
We can also obtain the same inequality for the other case with the same argument. Then, we get
\begin{align*}
\sup_{\substack{|v|=1\\|  w|=1}}F(x_{1},z_{1},v,w,g) \le M + \sup_{\substack{| v|=1\\|  w|=1}}F(x_{1},z_{1},v,w,f)
\end{align*}
and
\begin{align*}
\inf_{\substack{| v|=1\\| w|=1}}F(x_{1},z_{1},v,w,g) \le M + \inf_{\substack{|  v|=1\\|  w|=1}}F(x_{1},z_{1},v,w,f).
\end{align*}
Now from \eqref{add3}, we deduce that
\begin{align*}
u^{\eps}(x_{1})-u^{\eps}(z_{1})&=
\frac{1}{2}\sup_{\substack{| v|=1\\|  w|=1}}F(x_{1},z_{1},v,w,g)+\frac{1}{2}\inf_{\substack{|  v|=1\\| w|=1}}F(x_{1},z_{1}, v, w,g)
\\& \le M+ \frac{1}{2}\sup_{\substack{|  v|=1\\| w|=1}} F(x_{1},z_{1},v,w,f) + \frac{1}{2}\inf_{\substack{|  v|=1\\|  w|=1}}F(x_{1},z_{1}, v, w,f)
\\ & \le   u^{\eps}(x_{1})-u^{\eps}(z_{1})-f(x_{1},z_{1})+\eta
\\ & \qquad + \frac{1}{2}\sup_{\substack{| v|=1\\| w|=1}} F(x_{1},z_{1},v,w,f) + \frac{1}{2}\inf_{\substack{| v|=1\\|  w|=1}}F(x_{1},z_{1}, v, w,f),
\end{align*}
that is,
\begin{align*}
f(x_{1},z_{1}) \le  \frac{1}{2}\sup_{\substack{| v|=1\\| w|=1}} F(x_{1},z_{1},v,w,f) + \frac{1}{2}\inf_{\substack{|  v|=1\\| w|=1}}F(x_{1},z_{1}, v, w,f)+\eta.
\end{align*}
Thus, we can derive a contradiction if we show
\begin{align} \label{aim1}
f(x,z) >  \frac{1}{2}\sup_{\substack{|  v|=1\\| w|=1}} F(x,z,v,w,f) + \frac{1}{2}\inf_{\substack{| v|=1\\|  w|=1}}F(x,z, v, w,f)
\end{align}
for every $(x,z ) \in B_{1}\times B_{1}$.

The case $|x-z|\approx \eps$ follows from the fact that the steps are of size  $\eps$.
We include the details later.
Now we focus on the case $|x-z|>\frac{N}{10} \eps$.
In this case, since $f_{2}=0$, we need to prove that
\begin{align} 
\label{midest1}
f_{1}(x,z)>
\frac{1}{2}\sup_{\substack{|v|=1\\|w|=1}}
F(x,z,v,w,f_{1})
+
\frac{1}{2}\inf_{\substack{|\tilde v|=1\\|\tilde w|=1}}
F(x,z,\tilde v,\tilde w,f_{1}).
\end{align}

We will use the following Taylor's expansion for the function $f_1$,
\begin{equation}
\label{eq:taylor}
\begin{split}
f_1&(x+h_x,z+h_z)\\
&=f_1(x,z)+C\delta|x-z|^{\delta-1}(h_x-h_z)_V+2 \langle x+z, h_x+h_z\rangle \\
&\,\,+\frac{C}{2}\delta|x-z|^{\delta-2}\big((\delta-1)(h_x-h_z)^2_V+(h_x-h_z)^2_{V^{\perp}}\big)\\
&\,\,+|h_x+h_z|^2\,+\mathcal{E}_{x,z}(h_x,h_z),
\end{split}
\end{equation}
where $V$ is the space spanned by $x-z$, $(h_x-h_z)_V$  refers to the scalar projection onto $V$ i.e. $\frac{\langle h_x-h_z, x-z\rangle }{\abs{x-z}}$, and $(h_x-h_z)_{V^{\perp}}$ onto the orthogonal complement (see also \cite{luirop18,arroyohp17}). 

By Taylor's theorem, the error term satisfies
\begin{equation}
|\mathcal{E}_{x,z}(h_x,h_z)|\leq C|(h_x,h_z)|^3(|x-z|-2\eps)^{\delta-3}\,, 
\end{equation}
if $|x-z|> 2\eps\,$. Especially, if we choose 
\begin{equation}\label{Noletus}
N\ge \frac{100C}{\delta}\,,
\end{equation}
then in the case $\abs{x-z}> N\frac{\eps}{10} $  and $|h_x|,|h_z|\leq \eps$, we can see that $|x-z|-2\eps \ge \frac{|x-z|}{2}$ since
$$\abs{x-z}> \frac{10C\eps}{\delta} > 10\eps , $$
and thus we have

\begin{equation}
\label{error2}
\begin{split}
|\mathcal{E}_{x,z}(h_x,h_z)|&\leq C(2\eps)^3\bigg(\frac{|x-z|}{2}\bigg)^{\delta-3}\,
\\
&\leq 
64C\eps^2|x-z|^{\delta-2}\frac{\eps}{|x-z|}\\
&\leq 64\eps^2|x-z|^{\delta-2}\frac{\delta}{10}\\
&\leq 10\eps^2|x-z|^{\delta-2}.
\end{split}
\end{equation}
Observe that when considering the Taylor expansions, all the first order term will be canceled.

Now we are ready to proceed to the core of the matter, that is to prove \eqref{midest1}.
Here is where the precise definition of when to apply the rule (i) or (ii) plays the main role.
We will estimate the infimum of $F$ by considering
\[
\tilde v=\frac{x-z}{|x-z|} \quad \text{and} \quad \tilde w=-\frac{x-z}{|x-z|}.
\]
Observe that in this case, $v_{y^\perp}=w_{y^\perp}=0$, hence, the rule (i) applies.
Recall that
$$F(x,z,\tilde v,\tilde w,f) = \frac{ f((x,z) + \eps (\tilde v,\tilde w)) +f((x,z) - \eps (\tilde v,\tilde w)) }{2}$$
in this case.
From \eqref{eq:taylor}, we have
\begin{align*}
&F(x,z,\tilde v,\tilde w,f) -f(x,z) \\ &= \frac{C}{2}\delta|x-z|^{\delta-2}\big((\delta-1)(\tilde v-\tilde w)^2_y+(\tilde v-\tilde w)^2_{y^\perp}\big)+|\tilde v+\tilde w|^2+\mathcal{E}_{x,z}(\tilde v,\tilde w).
\end{align*}
Therefore, we obtain
\[
\begin{split}
&\inf_{\substack{|\tilde v|=1\\|\tilde w|=1}}
F(x,z,\tilde v,\tilde w,f) - f(x,z)\\
&\leq \frac{C}{2}\delta|x-z|^{\delta-2}\big((\delta-1)(\tilde v-\tilde w)^2_y+(\tilde v-\tilde w)^2_{y^\perp}\big)+|\tilde v+\tilde w|^2+\mathcal{E}_{x,z}(\tilde v,\tilde w)\\
&\leq \frac{C}{2}\delta|x-z|^{\delta-2}(\delta-1)(2\eps)^2+10\eps^2|x-z|^{\delta-2}.\\
\end{split}
\]
To bound the supremum, we have to separate the cases depending on whether rule (i) or rule (ii) applies. The key point here is to bound $(h_x-h_y)_{V^\perp}^2$ by strictly less than $(2\eps)^2$.

When the rule (ii) is applied, we have
\[
(v-(-w))_{y^\perp}^2=v_{y^\perp}^2+w_{y^\perp}^2+2\langle v_{y^\perp},w_{y^\perp}\rangle,
\]
and if $\langle v_{y^\perp},w_{y^\perp}\rangle<0$, then 
$(v-(-w))_{y^\perp}^2\leq 2$.
If rule (i) is applied and $\langle v_{y^\perp},w_{y^\perp}\rangle\geq 0$ the same calculation can be performed.
It remains to check the case where $|v_{y^{\perp}}|^2+|w_{y^{\perp}}|^2\leq 1$. In this case we have
\[
(v-w)_{y^\perp}^2=v_{y^\perp}^2+w_{y^\perp}^2-2\langle v_{y^\perp},w_{y^\perp}\rangle
\le 2(v_{y^\perp}^2+w_{y^\perp}^2) \le 2
\]
since $2| \langle v_{y^\perp},w_{y^\perp}\rangle| \le  v_{y^\perp}^2+w_{y^\perp}^2$.
Hence we can see that
\[
\begin{split}
&F(x,z,v,w,f) - f(x,z)\\
&\leq \frac{C}{2}\delta|x-z|^{\delta-2}\big((\delta-1)(v\pm w)^2_y+(v\pm w)^2_{y^\perp}\big)+|v\pm w|^2+\mathcal{E}_{x,z}(v,\pm w)\\
&\leq \frac{C}{2}\delta|x-z|^{\delta-2}2\eps^2+(2\eps)^2+10\eps^2|x-z|^{\delta-2}.\\
\end{split}
\]

Finally, we obtain
\begin{align*}
\begin{split}
&\sup_{\substack{|v|=1\\|w|=1}}
F(x,z,v,w,f)
+
\inf_{\substack{|\tilde v|=1\\|\tilde w|=1}}
F(x,z,\tilde v,\tilde w,f) - f(x,z)\\
&\leq \frac{C}{2}\delta|x-z|^{\delta-2}2\eps^2+(2\eps)^2+10\eps^2|x-z|^{\delta-2}
\\ & \qquad +\frac{C}{2}\delta|x-z|^{\delta-2}(\delta-1)(2\eps)^2+10\eps^2|x-z|^{\delta-2}\\
&\leq \frac{C}{2}\delta|x-z|^{\delta-2}((\delta-1)4\eps^2+2\eps^2)+4\eps^2+20\eps^2|x-z|^{\delta-2}\\
&\leq \eps^2|x-z|^{\delta-2}(2C\delta^{2}-C\delta+20)+4\eps^2.\\
\end{split}
\end{align*}

Now it is enough to show that
\begin{align}\label{est0}
\eps^2|x-z|^{\delta-2}(2C\delta^{2}-C\delta+20)+4\eps^2<0.
\end{align}
For $0<\delta<\frac{1}{2}$, we have 
\begin{align}
\label{eq:choice-of-delta}
\delta-2\delta^2>0.
\end{align}
Therefore, given $\tilde C>0$ we can take $C=\frac{\frac{\tilde C+4}{4^{\delta-2}}+20}{\delta-2\delta^2}>0$ so that
\[
2C\delta^{2}-C\delta+20
=-\frac{\tilde C+4}{4^{\delta-2}}.
\]
Observe that this is where we explicitly fix $\delta$.

Recalling that $|x-z|<4$, we obtain
\[
|x-z|^{\delta-2}\Big(-\frac{\tilde C+4}{4^{\delta-2}}\Big)+4 \leq -\tilde C.
\]

We have proved that
\begin{align} \label{wts1}
\frac{1}{2}\sup_{\substack{|v|=1\\|w|=1}}
F(x,z,v,w,f)
+
\frac{1}{2}\inf_{\substack{|\tilde v|=1\\|\tilde w|=1}}
F(x,z,\tilde v,\tilde w,f)-f(x,z)<-\tilde C\eps^2.
\end{align}

Now we consider the case $|x-z| \le \frac{N}{10} \eps$.
We remark that the counterassumption of \eqref{int2} cannot occur when $x=z$, since $$u(x)-u(x)-f(x,x)= C^{2N}\eps^{\delta}-4|x|^{2} \le  C^{2N}\eps^{\delta}.$$
Thus it is sufficient to show \eqref{aim1} when $0<|x-z| \le \frac{N}{10} \eps$.
We first observe that
\begin{align*}
\big||a|^{\delta}-|b|^{\delta}\big| \le |a-b|^{\delta}
\end{align*}
for any $a,b \in \mathbb{R}^{n}$ since $0<\delta<1 $. 
This yields
\begin{align*}
\big||(x+\eps v)-(z+ \eps w)|^{\delta}-|x-z|^{\delta}\big| \le \eps^{\delta}|v-w|^{\delta}.
\end{align*}
Therefore, we see that 
\begin{align*}
|f_{1}((x,z)+\eps (v,w))-f_{1}(x,z)| & \le 
C\eps^{\delta}|v-w|^{\delta}+ 2\eps |\langle x+z,v+w \rangle|+\eps^{2}|v+w|^{2}
\\ & \le 2C\eps^{\delta}+ 8\eps +4\eps^{2}
\\ & \le 3C\eps^{\delta}
\end{align*}
for any unit vectors $v$ and $w$, and sufficiently large $C$. 
Then we obtain
\begin{align*}
\sup_{\substack{|v|=1\\|w|=1}}F(x,z,v,w,f_{1})-f_{1}(x,z) \le 3C\eps^{\delta}
\end{align*}
in both cases.
Since $f_{2}\ge 0$ and $f=f_{1}-f_{2}$, we have
$$ \sup_{\substack{|v|=1\\|w|=1}}F(x,z,v,w,f) \le \sup_{\substack{|v|=1\\|w|=1}}F(x,z,v,w,f_{1}).$$
Observe that
$$  \inf_{\substack{|v|=1\\|w|=1}}F(x,z,v,w,f) \le  \sup_{\substack{|v|=1\\|w|=1}}F(x,z,v,w,f_{1})- \sup_{\substack{|v|=1\\|w|=1}}F(x,z,v,w,f_{2}).$$
Since, given $(x,z) \in A_{i}$, we can choose unit vectors $v,w $ such that $(x,z)+\eps(v,w) \in A_{i-1}$, we see that 
\begin{align*}
\sup_{\substack{|v|=1\\|w|=1}}F(x,z,v,w,f_{2})& \ge \frac{1}{2}\big( f_{2}((x,z)+\eps(v,w))\big)
\\ & = \frac{1}{2} C^{2(N-i+1)}\eps^{\delta}
\\ & = \frac{1}{2} C^{2(N-i)}\eps^{\delta}(C^{2}-4) + 2f_{2}(x,z).
\end{align*}
By choosing a large constant $C$, we get 
\begin{align} \label{supf2} \sup_{\substack{|v|=1\\|w|=1}}F(x,z,v,w,f_{2})  <  6C\eps^{\delta}+ 2f_{2}(x,z).
\end{align}
Combining the previous estimates, we see that
\begin{align*} \frac{1}{2}\sup_{\substack{|v|=1\\|w|=1}}
&F(x,z,v,w,f) 
+
\frac{1}{2}\inf_{\substack{|\tilde v|=1\\|\tilde w|=1}}
F(x,z,\tilde v,\tilde w,f) \\ &  <   \sup_{\substack{|v|=1\\|w|=1}}F(x,z,v,w,f_{1}) -
3C\eps^{\delta} - f_{2}(x,z)
\\ & \le (f_{1}(x,z)+3C\eps^{\delta})  -
3C\eps^{\delta} - f_{2}(x,z)
\\ & =f(x,z).
\end{align*}
This yields \eqref{aim1}.
\end{proof}

\section{The general case}
In this section, we consider the DPP \eqref{dpp2} related to the PDE \eqref{eq:mainN}.
We rewrite the equation and present the rules of the game in a slightly different way.
We assume $\alpha=2\min\{\alpha_1,\alpha_n\}>0$.
We define $\beta=1-\alpha$, $\beta_i=\alpha_i/\beta$ for $i=2,\dots,n-1$ and $\beta_i=(\alpha_i-\alpha/2)/\beta$ for $i=1,n$.
We can rewrite \eqref{eq:mainN} as
\begin{equation}
\label{rewritten}
\alpha\frac{\lambda_1+\lambda_n}{2}+\beta \sum_{i=1}^n \beta_i\lambda_i=0.
\end{equation}

We remark that one can derive H\"{o}lder regularity for \eqref{rewritten},
since a viscosity solution to \eqref{rewritten} satisfies
\begin{align*}
&\bigg(  1-\frac{\alpha}{2} \bigg)\sum_{\lambda_{i}<0}\lambda_{i}+\frac{\alpha}{2n}\sum_{\lambda_{i}>0}\lambda_{i} \le
\bigg(  1-\frac{(n-1)\alpha}{2n} \bigg)\lambda_{1}+\frac{\alpha}{2n} \sum_{i=2}^n \lambda_{i} \\ &\qquad\le \alpha\frac{\lambda_1+\lambda_n}{2}+\beta \sum_{i=1}^n \beta_i\lambda_i \\ &  \le 
\bigg(  1-\frac{(n-1)\alpha}{2n} \bigg)\lambda_{n}+\frac{\alpha}{2n} \sum_{i=1}^{n-1} \lambda_{i} \le \bigg(  1-\frac{\alpha}{2} \bigg)\sum_{\lambda_{i}>0}\lambda_{i}+\frac{\alpha}{2n}\sum_{\lambda_{i}<0}\lambda_{i} .
\end{align*} 
These Pucci type inequalities are what is required to use \cite[Proposition 4.10]{caffarelli95book}.

Now the game for $\alpha_i$ can be presented in the following way: at every round with probability $\alpha$ we play the game for $\frac{1}{2}\lambda_1+\frac{1}{2}\lambda_n$ and with probability $\beta$ we play the game according to $\beta_i$.
In this case, the related DPP is
\begin{align}
\label{dpp2var}\begin{split}
u^\eps (x) =&
\frac{\alpha}{2}\sup_{|v|=1} \frac{ u^\eps (x + \eps v) + u^\eps (x - \eps v)}{2}
+\frac{\alpha}{2}\inf_{|w|=1} \frac{ u^\eps (x + \eps w) + u^\eps (x - \eps w)}{2}
\\ & \quad  + \beta  \sum_{i=1}^n \beta_i \inf_{dim(S)=j}\sup_{\substack{v\in S\\ |v|=1}}
\frac{ u^\eps (x + \eps v) + u^\eps (x - \eps v)}{2}.
\end{split}
\end{align}
This is equivalent to \eqref{dpp2}.

In order to define the $2n$-dimensional game related to \eqref{dpp2var}, first we define a $2n$-dimensional game related to $\lambda_j$.
Fix $j\in \{1,...,n\}$.
We consider the game related to $\lambda_{j}$ and write $u_{j}$ for its value function. It satisfies the following DPP
\begin{align*}
u_{j}(x) = \inf_{dim(S)=j}\sup_{\substack{v\in S\\ |v|=1}}
\frac{ u_{j} (x + \eps v) +u_{j} (x - \eps v)}{2}
\end{align*}
for any $x \in \Omega$.

Set $g_{j}(x,z)=u_{j}(x)-u_{j}(z)$.
We have
\begin{align*}
&\inf_{dim(S)=j}\sup_{\substack{v\in S\\ |v|=1}}
\frac{ u_{j} (x + \eps v) +u_{j} (x - \eps v)}{2} 
- \inf_{dim(\tilde S)=j}\sup_{\substack{\tilde v\in \tilde S\\ |\tilde v|=1}}
\frac{ u_{j}(z + \eps \tilde v) +u_{j} (z - \eps \tilde v)}{2}
\\ 
& = \hspace{-0.6em} \sup_{dim(\tilde S)=j}\inf_{dim(S)=j}\sup_{\substack{v\in S\\ |v|=1}}\inf_{\substack{\tilde v\in \tilde S\\ |\tilde v|=1}} \bigg\{
  \frac{ u_{j} (x + \eps v) +u_{j} (x - \eps v)}{2} 
-
\frac{ u_{j}(z + \eps \tilde v) +u_{j} (z - \eps \tilde v)}{2}\bigg\}
\\ 
&=  
\sup_{dim(\tilde S)=j}\inf_{dim(S)=j}\sup_{\substack{v\in S\\ |v|=1}}\inf_{\substack{\tilde v\in \tilde S\\ |\tilde v|=1}}
\frac{g_{j}(x + \eps v,z + \eps \tilde v)+g_{j}(x - \eps v,z - \eps \tilde v)}{2}.
\end{align*}
We can read the rules of the $2n$-dimensional game as follows:
 Player~II selects $\tilde S$, Player~I the subspace $S$ and then Player~II a unitary vector $v\in S$ and Player~I a vector $\tilde v\in \tilde S$.
Then, the token moves to $(x,z)\pm\eps (v,\tilde v)$ each with probability one half.

Combining the above observation for each $g_{j}$ with the $2n$-dimensional DPP for the game associated with $\frac{1}{2}\lambda_1+\frac{1}{2}\lambda_n$, for the function $g(x,z)=u^{\eps}(x)-u^{\eps}(z)$, we have
\begin{align*}
g(x,z)& =\frac{\alpha}{2}\sup_{\substack{|v|=1\\|w|=1}}F(x,z,v,w,g)
 +
\frac{\alpha}{2}\inf_{\substack{|\tilde v|=1\\|\tilde w|=1}}F(x,z,v,w,g)
\\ & + \beta  \sum_{i=1}^n \beta_i \sup_{dim(\tilde S)=j}\inf_{dim(S)=j}\sup_{\substack{v\in S\\ |v|=1}}\inf_{\substack{\tilde v\in \tilde S\\ |\tilde v|=1}}
\frac{g(x + \eps v,z + \eps \tilde v)+g(x - \eps v,z - \eps \tilde v)}{2},
\end{align*}
where $F$ is the function given by \eqref{defcapf}.

Now we state and prove the H\"{o}lder regularity result for \eqref{dpp2}.
\begin{theorem}
\label{thm:comb-holder}
Let $u^{\eps}$ be a function satisfying the DPP \eqref{dpp2} in a bounded domain $\Omega$.
Then for any $0<\delta<\frac{1}{2}$ and $x,z\in B_{r}$ with $B_{2r} \subset\Omega$, there exists a constant $C=C(\delta, \alpha)>0$ such that
$$|u^{\eps}(x)-u^{\eps}(z)| \le C||u^{\eps}||_{L^{\infty}(B_{2r})}\bigg(\frac{|x-z|^{\delta}}{r^{\delta}}+\frac{\eps^{\delta}}{r^{\delta}}\bigg). $$ 
\end{theorem}

\begin{proof}
Recall the barrier function $f$ in the proof of Theorem \ref{thm:main}.
By a similar argument as in the previous section, it is enough to show that
\begin{align} \label{dpp3}\begin{split}
&f(x,z) >\frac{\alpha}{2}\sup_{\substack{|v|=1\\|w|=1}}F(x,z,v,w,f)  +
\frac{\alpha}{2}\inf_{\substack{|v|=1\\| w|=1}}F(x,z,v,w,f)
\\ & \quad + \beta  \sum_{i=1}^n \beta_i \sup_{dim(\tilde S)=j}\inf_{dim(S)=j}\sup_{\substack{v\in S\\ |v|=1}}\inf_{\substack{\tilde v\in \tilde S\\ |\tilde v|=1}}
\frac{f(x + \eps v,z + \eps \tilde v)+f(x - \eps v,z - \eps \tilde v)}{2}.
\end{split}\end{align}

We first consider the case $|x-z|>\frac{N}{10} \eps$.
For the terms involved in the game associated with $\frac{1}{2}\lambda_1+\frac{1}{2}\lambda_n$ we can recall the estimate \eqref{wts1}.
Meanwhile, for the game associated to $\lambda_j$, 
 we observe that by taking $S=\tilde S$ and $\tilde v=v$ we get
\begin{align}  \label{aa1}
\begin{split}
\sup_{dim(\tilde S)=j}\inf_{dim(S)=j}\sup_{\substack{v\in S\\ |v|=1}}\inf_{\substack{\tilde v\in \tilde S\\ |\tilde v|=1}}
\frac{f(x + \eps v,z + \eps \tilde v)+f(x - \eps v,z - \eps \tilde v)}{2}\\ \le
\sup_{dim(\tilde S)=j}\sup_{\substack{\tilde{v}\in S\\ |\tilde{v}|=1}}
\frac{f(x + \eps \tilde v,z + \eps \tilde v)+f(x - \eps \tilde v,z - \eps \tilde v)}{2}.
\end{split}
\end{align}
Moreover, we observe that
\begin{align} \label{aa2}
\sup_{dim(\tilde S)=j}\sup_{\substack{\tilde{v}\in S\\ |\tilde{v}|=1}}\frac{f(x + \eps \tilde v,z + \eps \tilde v)+f(x - \eps \tilde v,z - \eps \tilde v)}{2}
=
f(x,z)+4\eps^2.
\end{align}
We conclude that
\begin{align*}
&\frac{\alpha}{2}\sup_{\substack{|v|=1\\|w|=1}}F(x,z,v,w,f)  +
\frac{\alpha}{2}\inf_{\substack{|\tilde v|=1\\|\tilde w|=1}}F(x,z,v,w,f)
\\ & \quad + \beta  \sum_{i=1}^n \beta_i \sup_{dim(\tilde S)=j}\inf_{dim(S)=j}\sup_{\substack{v\in S\\ |v|=1}}\inf_{\substack{\tilde v\in \tilde S\\ |\tilde v|=1}}
\frac{f(x + \eps v,z + \eps \tilde v)+f(x - \eps v,z - \eps \tilde v)}{2}
\\ & <  f(x,z) -\alpha \tilde C\eps^2 + 4 \beta\eps^2,
\end{align*}
where $\tilde{C}$ is the constant in $\eqref{wts1}$.
Thus, if we take $\tilde C$ large enough such that
\[
-\tilde C\alpha+4\beta<0,
\]
we obtain \eqref{dpp3}.

Next we assume that $|x-z|\le \frac{N}{10} \eps$.
From \eqref{dpp3}, \eqref{aa1} and \eqref{aa2}, it is enough to show that 
\begin{align*} 
f(x,z) >\frac{\alpha}{2}\sup_{\substack{|v|=1\\|w|=1}}F(x,z,v,w,f)  +
\frac{\alpha}{2}\inf_{\substack{| v|=1\\| w|=1}}F(x,z,v,w,f)
+ \beta \sum_{i=1}^n \beta_i( f(x,z)+4\eps^2).
\end{align*}
This can be rewritten as
\begin{align} \label{aa3}
f(x,z) >\frac{1}{2}\sup_{\substack{|v|=1\\|w|=1}}F(x,z,v,w,f)  +
\frac{1}{2}\inf_{\substack{| v|=1\\| w|=1}}F(x,z,v,w,f)
+\frac{4(1-\alpha)}{\alpha}\eps^2.
\end{align}
We use a similar argument in the proof of Theorem \ref{thm:main}, but we choose $C$ sufficiently large such that
\begin{align*} \sup_{\substack{|v|=1\\|w|=1}}F(x,z,v,w,f_{2}) \ge \bigg(\frac{8(1-\alpha)}{\alpha} + 6 \bigg) C\eps^{\delta}+ 2f_{2}(x,z)
\end{align*}
in \eqref{supf2}.
Then we get \eqref{aa3}, which completes the proof.
\end{proof}

\subsection{The Dominative $p$-Laplacian}

Recently, there has been some interest to the Dominative $p$-Laplacian
$$\mathcal{D}_{p}u:=\lambda_{1}+ \cdots+ \lambda_{n-1}+ (p-1)\lambda_{n}, $$
where $2 \le p < \infty$, see \cite{brustad2020superposition}.
It explains a superposition of $p$-superharmonic functions, which was studied in \cite{MR2350398}.

The dominative $p$-Laplacian can be regarded as a special case of the operator $\sum_{i=1}^n \alpha_i\lambda_i$,
which has been considered so far.
Observe that the equation $\mathcal{D}_{p}u=0$ is equivalent to the equation \eqref{eq:mainN} when
$\alpha_{i}=\frac{1}{n+p-2}$ for $i=1,\cdots,n-1$, and $\alpha_{n}=\frac{p-1}{n+p-2}$.
Therefore, by plugging these values in $\eqref{dpp2}$ we obtain the following DPP
\begin{align*}
u^\eps (x) =&\frac{1}{n+p-2}\sum_{j=1}^{n-1} \inf_{dim(S)=j}\sup_{\substack{v\in S\\|v|=1}}
\frac{ u^\eps (x + \eps v) + u^\eps (x - \eps v)}{2}\\
&\quad +\frac{p-1}{n+p-2}\sup_{\substack{v\in S\\|v|=1}}
\frac{ u^\eps (x + \eps v) + u^\eps (x - \eps v)}{2}.
\end{align*}
Since $\min\{\alpha_1,\alpha_n\}>0$, our result, Theorem \ref{thm:domp}, covers the solutions to this DPP.

We also remark here that the operator $\mathcal{D}_{p}$ is uniformly elliptic, and thus we can obtain $C^{1,\delta}$-regularity for viscosity solutions of the equation $\mathcal{D}_{p}u=0$ (see \cite[Section 5.3]{caffarelli95book}).

A different game associated to the Dominative $p$-Laplacian was presented in \cite{brustadlm20} (see also \cite{hoeg2020control}).
Their DPP reads as follows
\begin{align} \label{dppdomp}
u^{\eps}(x)= q \kint_{B_{\eps}(x)}u^{\eps}(y) dy +  (1-q)\sup_{|v|=1}\left\{ \frac{ u^\eps (x + \eps v) + u^\eps (x - \eps v)}{2}\right\},
\end{align}
where $q=\frac{n+2}{n+p}$.
This is a control problem.
Let $x_{0}$ be the starting point.
The player first chooses a unit vector $v$, and then the token is moved according to the following rules:
$x_{1}$ is randomly selected in $B_{\eps}(x_{0})$ with probability $\frac{n+2}{n+p}$, and
$x_{1}=x_{0}\pm \eps v$ with probability $\frac{p-2}{2(n+p)}$, respectively.
This stochastic process is repeated until the token leaves $\Omega$. 
The player tries to maximize the expected value of $G(x_{\tau})$, and thus he/she chooses the direction $v$ for this purpose.

We can also obtain the following regularity result for this DPP with a slightly worse upper bound for $\delta$.
\begin{theorem}
\label{thm:domp}
Let $u^{\eps}$ be a function satisfying the DPP \eqref{dppdomp} in a bounded domain $\Omega$.
Then for any $0<\delta<\frac{1}{10}$ and $x,z\in B_{r}$ with $B_{2r} \subset\Omega$, there exists a constant $C=C(\delta)>0$ such that
$$|u^{\eps}(x)-u^{\eps}(z)| \le C||u^{\eps} ||_{L^{\infty}(B_{2r})}\bigg(\frac{|x-z|^{\delta}}{r^{\delta}}+\frac{\eps^{\delta}}{r^{\delta}}\bigg). $$ 
\end{theorem}
\begin{proof}
We first observe that
\begin{align*}
&\sup_{|v|=1}\left\{ \frac{ u^\eps (x + \eps v) + u^\eps (x - \eps v)}{2}\right\}-\sup_{|w|=1}\left\{ \frac{ u^\eps (z + \eps w) + u^\eps (z - \eps w)}{2}\right\}
\\ & = \sup_{|v|=1}\inf_{|w|=1}\left\{ \frac{ u^\eps (x + \eps v) + u^\eps (x - \eps v)- u^\eps (z + \eps w) - u^\eps (z - \eps w)}{2}\right\}.
\end{align*}
Like before, we again  consider $$g(x,z)=u^{\eps}(x)-u^{\eps}(z) .$$
Then we have
\begin{align*}
&g(x,z)\\&=u^{\eps}(x)-u^{\eps}(z)
\\ & = q \bigg( \kint_{B_{\eps}(x)}u^{\eps}( y)dy -\kint_{B_{\eps}(z)}u^{\eps}( y)dy \bigg)
\\ & \quad +(1-q) \sup_{|v|=1}\inf_{|w|=1}\left\{ \frac{ u^\eps (x + \eps v) + u^\eps (x - \eps v)- u^\eps (z + \eps w) - u^\eps (z - \eps w)}{2}\right\}
\\& =q \bigg( \kint_{B_{\eps}(x)}u^{\eps}( y)dy -\kint_{B_{\eps}(z)}u^{\eps}( y)dy \bigg)
\\ & \qquad +(1-q)  \sup_{|v|=1}\inf_{|w|=1}\left\{ \frac{g ((x,z)+ \eps( v,w)) + g ((x,z)- \eps( v,w)) }{2}\right\}
\\ & \le q  \bigg( \kint_{B_{\eps}(x)}u^{\eps}( y)dy -\kint_{B_{\eps}(z)}u^{\eps}( y)dy \bigg)
\\ & \qquad+(1-q) \sup_{|v|=1}\left\{ \frac{g ((x,z)+ \eps( v,v)) + g ((x,z)- \eps( v,v)) }{2}\right\}.
\end{align*}

Again, we recall the auxiliary function $f$ and the ideas explained before.
In \cite[Section 4]{luirop18}, we can find the following observation 
\begin{align*} &  \kint_{B_{\eps}(x)}u^{\eps}( y)dy -\kint_{B_{\eps}(z)}u^{\eps}( y)dy
\\ & = \frac{1}{|B_{\eps}|}\bigg( \int_{B_{\eps}(0)\backslash B_{\eps}(z-x)}\big(u(x+h)-u(z+P_{x,z}(h)) \big) dh
 \bigg), \end{align*}
where 
 $P_{x,z}$ is a map to send a point to its mirror point with respect to $\textrm{span}(x-z)^{\perp}$, the orthogonal complement of the subspace generated by $x-z$.
Repeating a similar calculation in the proof of Theorem \ref{thm:main}, we see that
it is enough to deduce a contradiction to \eqref{as3} if we prove
\begin{align} \label{add4} \begin{split}
f(x,z)& >  q  \cdot \frac{1}{|B_{\eps}|}\bigg( \int_{B_{\eps}(0)\backslash B_{\eps}(z-x)}f(x+h,z+P_{x,z}(h)) dh \\& \qquad \qquad \qquad
 + \int_{B_{\eps}(x)\cap B_{\eps}(z)}f(y,y) dy \bigg) 
\\ & + (1-q) \sup_{|v|=1}\left\{ \frac{f((x,z)+ \eps( v,v)) +f ((x,z)- \eps( v,v)) }{2}\right\} 
\end{split}\end{align}
for every $(x.z) \in B_{1}\times B_{1}$ (see also (4.25) in \cite{luirop18}).

Assume $\delta <1/10$ and set $ C=\frac{10^{10}}{\delta^{2}\omega}$, where $\omega$ is to be determined.
We refer to the following estimate in \cite[Section 4]{luirop18}: 
\begin{align*} 
&\frac{1}{|B_{\eps}|}\bigg( \int_{B_{\eps}(0)\backslash B_{\eps}(z-x)} \hspace{-1em}f(x+h,z+P_{x,z}(h)) dh
 + \int_{B_{\eps}(x)\cap B_{\eps}(z)}\hspace{-1em} f(y,y) dy \bigg) -f(x,z)
\\ & < K \eps^{\delta},
\end{align*}
where
\begin{align*}K= \left\{ \begin{array}{ll}
|x-z|^{\delta-2}\big( 10-\frac{C\delta}{4(n+2)}\big)& \textrm{if $|x-z| >  N\eps/10, $}\\
\big( -\frac{C^{2}}{4^{n}}+3C+1 \big) & \textrm{if $|x-z| \le N\eps / 10$.} 
\end{array} \right.
\end{align*}
Thus, by choosing $\omega \le 4^{-n}$, we get
\begin{align} \label{lpest} \begin{split}
&\frac{1}{|B_{\eps}|}\bigg( \int_{B_{\eps}(0)\backslash B_{\eps}(z-x)} \hspace{-1em}f(x+h,z+P_{x,z}(h)) dh
 + \int_{B_{\eps}(x)\cap B_{\eps}(z)}\hspace{-1em} f(y,y) dy \bigg) -f(x,z)
\\ & < -\tilde{C} \eps^{\delta}
\end{split}\end{align}
for $\tilde{C}=\min \{ \frac{1}{4}( \frac{C\delta}{4(n+2)}-10),  \frac{C^{2}}{4^{n}}-3C-1\}$.
Nothe that $\tilde{C}$ is strictly positive because  $ C=\frac{10^{10}}{\delta^{2}\omega}$.
We also observe that 
\begin{align*}
\sup_{|v|=1}\left\{ \frac{f((x,z)+ \eps( v,v)) +f ((x,z)- \eps( v,v)) }{2}\right\} =f(x,z)+4\eps^2
\end{align*}
for any $(x,z) \in B_{1}\times B_{1}$.
Combining this with \eqref{lpest}, we have
\begin{align*}
&   q  \cdot \frac{1}{|B_{\eps}|}\bigg( \int_{B_{\eps}(0)\backslash B_{\eps}(z-x)}f(x+h,z+P_{x,z}(h)) dh
 + \int_{B_{\eps}(x)\cap B_{\eps}(z)}f(y,y) dy \bigg) 
\\ & + (1-q) \sup_{|v|=1}\left\{ \frac{f((x,z)+ \eps( v,v)) +f ((x,z)- \eps( v,v)) }{2}\right\} -f(x,z)
\\ & < -q\tilde{C} \eps^{\delta} +4(1-q)\eps^2
\\& < \big(-q\tilde{C} +4(1-q)\big)\eps^{\delta}.
\end{align*}
Now we can complete the proof if we choose $\omega $ small such that $$-q\tilde{C} +4(1-q)<0 .$$
Combining the above estimates, we obtain \eqref{add4}. 
\end{proof}

\medskip

{\bf Acknowledgements.} 
The authors would like to thank Julio D. Rossi for useful discussions.


{\bf Funding.}
J.~H.\ was supported by NRF-2021R1A6A3A14045195. P.~B.\ and M.~P.\ were partly supported by the Academy of Finland project 298641. 


{\bf Declaration of Competing Interests.}
The authors declare that there is no conflict of interest regarding the publication of this article.


{\bf Data availability.}
This manuscript has no associated data.
~
\bibliographystyle{alpha}

\newcommand{\etalchar}[1]{$^{#1}$}
\def\cprime{$'$} \def\cprime{$'$} \def\cprime{$'$}

\end{document}